\documentclass[12pt]{amsart}
\usepackage{amsmath, amsthm, amsfonts, amssymb, color}
 \usepackage{mathrsfs}
\usepackage{amsfonts, amsmath,amsrefs}
 \usepackage{amsmath,amstext,amsthm,amssymb,amsxtra}
 \usepackage{txfonts} %also pxfonts
 \usepackage[colorlinks, citecolor=blue,hypertexnames=false]{hyperref}
 \allowdisplaybreaks
 \usepackage{pgf}
 \usepackage{multirow}%newtablepackage1
 \usepackage{diagbox} %newtablepackage2

 \usepackage{tikz} 
\usetikzlibrary{decorations.pathreplacing}

\allowdisplaybreaks

 \usepackage{comment}

\usepackage{mathtools}
\mathtoolsset{showonlyrefs,showmanualtags}

 \DeclareFontFamily{U}{wncy}{}
\DeclareFontShape{U}{wncy}{m}{n}{<->wncyr10}{}
\DeclareSymbolFont{mcy}{U}{wncy}{m}{n}
\DeclareMathSymbol{\Sh}{\mathord}{mcy}{"58}

\begin{document}

 \baselineskip 16.6pt
\hfuzz=6pt

\widowpenalty=10000

\newtheorem{cl}{Claim}
\newtheorem{theorem}{Theorem}[section]
\newtheorem{proposition}[theorem]{Proposition}
\newtheorem{coro}[theorem]{Corollary}
\newtheorem{lemma}[theorem]{Lemma}
\newtheorem{definition}[theorem]{Definition}
\newtheorem{assum}{Assumption}[section]
\newtheorem{example}[theorem]{Example}
\newtheorem{remark}[theorem]{Remark}
\renewcommand{\theequation}
{\thesection.\arabic{equation}}

\def\SL{\sqrt H}

\newcommand{\cent}{\operatorname{cent}}

\newcommand{\wid}{\operatorname{width}}
\newcommand{\heit}{\operatorname{height}}

\newcommand{\cdim}{n}

\newcommand{\set}[1]{\mathfrak{#1}}
\newcommand{\vrect}{\set{V}}
\newcommand{\Stack}{\set{S}}
\newcommand{\tile}{\set{T}}

\newcommand{\mar}[1]{{\marginpar{\sffamily{\scriptsize
        #1}}}}

\newcommand{\as}[1]{{\mar{AS:#1}}}
\newcommand\C{\mathbb{C}}
\newcommand\Z{\mathbb{Z}}
\newcommand\R{\mathbb{R}}
\newcommand\RR{\mathbb{R}}
\newcommand\CC{\mathbb{C}}
\newcommand\NN{\mathbb{N}}
\newcommand\ZZ{\mathbb{Z}}
\newcommand\HH{\mathbb{H}}
\def\RN {\mathbb{R}^n}
\renewcommand\Re{\operatorname{Re}}
\renewcommand\Im{\operatorname{Im}}

\newcommand{\mc}{\mathcal}
\newcommand\D{\mathcal{D}}
\def\hs{\hspace{0.33cm}}
\newcommand{\la}{\alpha}
\def \l {\alpha}
\newcommand{\eps}{\tau}
\newcommand{\pl}{\partial}
\newcommand{\supp}{{\rm supp}{\hspace{.05cm}}}
\newcommand{\x}{\times}
\newcommand{\lag}{\langle}
\newcommand{\rag}{\rangle}

\newcommand{\lset}{\left\lbrace}
\newcommand{\rset}{\right\rbrace}

\newcommand\wrt{\,{\rm d}}

\title[Schatten class in the Bloom setting]{SCHATTEN CLASSES AND COMMUTATOR IN THE TWO WEIGHT SETTING,\\
 I.  HILBERT TRANSFORM}

\author{Michael Lacey}
\address{Michael Lacey, Department of Mathematics, Georgia Institute of Technology
Atlanta, GA 30332, USA}
\email{lacey@math.gatech.edu}

\author{Ji Li}
\address{Ji Li, School of Mathematical and Physical Sciences, Macquarie University, NSW 2109, Australia}
\email{ji.li@mq.edu.au}

\author{Brett D. Wick}
\address{Brett D. Wick, Department of Mathematics\\
         Washington University - St. Louis\\
         St. Louis, MO 63130-4899 USA
         }
\email{wick@math.wustl.edu}

  \date{}


\begin{abstract}
We characterize the Hilbert--Schmidt class membership of commutator with the Hilbert transform in the two weight setting. 
The characterization depends upon the symbol of the commutator being in a  new weighted  Besov space. 
 This follows from a Schatten class $S_p$ result for dyadic paraproducts, where $1< p < \infty  $.  
We discuss the difficulties in extending the dyadic result to the full range of Schatten classes for the Hilbert transform.   
\end{abstract}

\maketitle
\tableofcontents 

\section{Introduction and Statement of Main Results}
\setcounter{equation}{0}

Consider the commutator with Hilbert transform $H$ defined as follows.
$$[b,H](f)(x):= b(x)H(f)(x) - H(bf)(x).  $$
In this paper, we are focused on the two weight setting, given as follows.  

Suppose $1<p<\infty$, $\mu,\lambda$ in the Muckenhoupt weight class $A_p$, and $\nu= \lambda^{-{1\over p}}\mu^{1\over p}$.
 Bloom \cite{MR805955} showed that $[b,H]$ is bounded from $L^p_\mu(\mathbb R)$ to $L^p_\lambda(\mathbb R)$
 if and only if $b$ is in ${\rm BMO}_\nu(\mathbb R)$. 
 Just recently, the first and second authors \cite{MR4345997} proved that $[b,H]$ is compact from $L^p_\mu(\mathbb R)$ to $L^p_\lambda(\mathbb R)$
 if and only if $b$ is  in ${\rm VMO}_\nu(\mathbb R)$.

We characterize the Hilbert--Schmidt property of the commutator above in the two weight setting.  Key to this is a  the Schatten class membership of the dyadic variants of the commutator. 

The two weight case introduces   difficulties not present in the unweighted case.  Our characterizations depend upon new weighted Besov spaces. 
For the dyadic case, we find a convenient characterization for the dyadic commutators to be in any Schatten class for $1<p< \infty $. 
Those classes $S^p$ are defined 
 in \S \ref{s:2}. See \eqref{e:SpNorm}. 
 But for the continuous case, we have to content ourselves with just the case $p=2$, 
 that is the Hilbert-Schmidt case.  
 We state our results here. 
Membership of commutators in Schatten classes is characterized by the symbol of the commutator being in Besov spaces. We introduce  a new type of weighted Besov spaces with respect to the two weights $\mu$ and $\lambda$ as follows.

\begin{definition}
Suppose $\mu,\lambda\in A_2$ and set $\nu= \mu^{1\over2}\lambda^{-{1\over2}}$. Let $b\in L_{\rm loc}^{1}(\mathbb R)$. Then we say that $b$ belongs to the weight  Besov space $B_{\nu}^{2}(\mathbb R)$ if
$$ \|b\|_{B_{\nu}^{2}(\mathbb R)}=  \Bigg(\int_{ \mathbb R} \int_{\mathbb R}{\big|b(x)-b(y)\big|^2\over |x-y|^2}\lambda(x)\mu^{-1}(y)dydx \Bigg)^{1\over 2}<\infty.  $$
\end{definition}

\begin{theorem}\label{schatten H}
Suppose $\mu,\lambda\in A_2$ and set $\nu= \mu^{1\over2}\lambda^{-{1\over2}}$.
Suppose $b\in L_{\rm loc}^{1}(\mathbb R)$. Then commutator $[b,H]$ belongs to $S^2(L^2_\mu(\mathbb R),L^2_\lambda(\mathbb R))$ if and only if
$b\in B_{\nu}^{2}(\mathbb R)$. Moreover,  in this case we have
$$
\|b\|_{B_{\nu}^{2}(\mathbb R)}\approx \|[b,H]\|_{S^2(L^2_\mu(\mathbb R),L^2_\lambda(\mathbb R))}.
$$
\end{theorem}

The key new tools to prove  this result are 
\begin{enumerate}
  \item For sufficiency, the characterisation of dyadic paraproducts $\Pi_b$ as we establish in Theorem \ref{thm Pi and Haar multiplier}.  This characterization is 
  in terms of the dyadic Besov spaces, which contain the continuous version.  

\item The intersection of shifted dyadic weighted Besov spaces is the continuous weighted Besov space, mimicking a well known property of BMO space.   This addresses sufficiency in the continuous case.

\item   For necessity in the continuous case,  the equivalent characterisations of the dyadic weighted Besov norm and the decomposition of dyadic interval via median value.

\end{enumerate}
We emphasize that this strategy also depends upon the Hilbert space structure. The `obvious' generalization to $L^p$ is problematic. We return to this point in the concluding section of the paper.

We turn to the dyadic setting.  
Throughout, $\mathcal D$ is a choice of a dyadic system or grid in $\mathbb{R}$. 
Let $\{h_I \colon I\in \mathcal D\}$ be the Haar functions associated to $ \mathcal D$. 
For function, or symbol, $b$, let $\widehat b (I) = \langle b, h_I\rangle$ be the 
Haar coefficients of $b$.   The Besov spaces are then given by this definition.

\begin{definition}\label{def dy Bp}
{Suppose $\nu\in A_2$, $0<p< \infty$}. Let $b\in L_{\rm loc}^{1}(\mathbb R)$. Let $\mathcal{D}$ be an arbitrary dyadic interval system in $\mathbb R$. Then we say that $b$ belongs to the weight dyadic  Besov space $B_{\nu,d}^{p}(\mathbb R, \mathcal D) = B_{\nu,d}^{p}(\mathbb R)$ if
\begin{align*}
{ \|b\|_{B_{\nu,d}^{p}(\mathbb R)}:=\Bigg(\sum_{I\in\mathcal D} \bigg( { |\widehat{b}(I)|\, |I|^{1\over2}\over {\nu(I)}} \bigg)^p\Bigg)^{1\over p} <\infty.}
\end{align*}
\end{definition}

The paraproduct operator with symbol  $b$ is defined to be
\begin{equation}
\label{e:Paraproduct}
\Pi_b\equiv \sum_{I\in\mathcal{D}} \widehat{b}(I) h_I\otimes \frac{\mathsf{1}_I}{\left\vert I\right\vert}. 
\end{equation}
A Haar multiplier is defined  in terms of an arbitrary function $ \epsilon \colon \mathcal D \to \{ -1, 1\} $.  
\begin{equation} \label{e.Multiplier}
T_\epsilon=  \sum _{I\in \mathcal D} \epsilon_I \widehat f (I) h_I. 
\end{equation}
The main result in the dyadic case holds for all Schatten classes $1< p < \infty $. 

\begin{theorem}\label{thm Pi and Haar multiplier}
Suppose $1<p<\infty$, $\mu,\lambda\in A_2$ and set $\nu= \mu^{1\over2}\lambda^{-{1\over2}}$.
 Suppose   $b\in L^1_{\rm loc} (\mathbb R)$. 
 
{\rm(1)} $\Pi_b$ is in $S^p(L^2_\mu(\R) , L^2_\lambda(\R))$ if and only if $b$ is in $B_{\nu,d}^p(\R)$.  Moreover,
$$\|\Pi_b\|_{ S^p(L^2_\mu(\R) , L^2_\lambda(\R)) } \approx \|b\|_{B_{\nu,d}^p(\R)}.$$

{\rm (2)} If $b$ is in $B_{\nu,d}^p(\R)$ then $[b,T_\epsilon]$ is in $S^p(L^2_\mu(\R) , L^2_\lambda(\R))$. Conversely, suppose the coefficients $\{\epsilon_I\}_{I\in \mathcal D}$ satisfies the following non-degenerate condition: for each dyadic interval $Q\in \mathcal D$, there is another dyadic interval $\hat Q\in \mathcal D$ with $|Q|=|\hat Q|$ and ${\rm dist}(Q,\hat Q)\leq |Q|$, such that there exists $I_0\in\mathcal D$ satisfying $Q,\hat Q\subset I_0$,  $|I_0|\leq C|Q|$, and
$$ \bigg|\sum_{I\in\mathcal D, I_0\subset I} \epsilon_I h_I(x)h_I(y)\bigg|\geq {C_0 \over |I_0|}, $$
where $C$ and $C_0$ are absolute positive constants.  Then, if $[b,T_\epsilon]$ is in $S^p(L^2_\mu(\R) , L^2_\lambda(\R))$, we deduce that $b$ is in $B_{\nu,d}^p(\R)$.
\end{theorem}

\smallskip
We now establish the dyadic structure for the weighted Besov spaces, in the case $p=2$. 
\begin{theorem}\label{thm dyadic structure}
Suppose $\mu,\lambda\in A_2$ and set $\nu= \mu^{1\over2}\lambda^{-{1\over2}}$.  There are two dyadic systems $ \mathcal D ^{0}$ and $ \mathcal D ^{1}$ for which we have 
 $ B_{\nu}^{2}(\mathbb R) =B_{\nu,d}^{2}(\mathbb R, \mathcal D ^{0})\cap B_{\nu,d}^{2}(\mathbb R, \mathcal D ^{1}) $ with 
$$ \|b\|_{B_{\nu}^{2}(\mathbb R)} \approx \|b\|_{B_{\nu,d}^{2}(\mathbb R ,\mathcal D ^{0})}+\|b\|_{B_{\nu,d}^{2}(\mathbb R, \mathcal D ^{1})}. $$
 
\end{theorem}

There is no mystery as to the selection of the two dyadic grids. The `one-third' shift suffices. See for example \cite{MR3420475}.

\bigskip 

In the unweighted case, these results are classical in nature.  
 The case of the Hilbert transform commutator is equivalent to Schatten class membership for Hankel operators.  This was characterized in terms of membership of the symbol in a Besov space  by  Peller \cites{MR4111756,MR1949210}.  
 The characterization holds for all $0< p < \infty$.  
The higher dimensional case $\R^n$, $n\geq2$, was studied by Janson and Wolff \cite{MR686178}.  
Interestingly, for the Riesz transform commutators, one has  membership in $S^p$ for 
$0< p \leq n$ if and only if the symbol is constant.  
Due to our difficulties in the continuous case, 
it is natural for this paper to be restricted to dimension one, and the Hilbert transform kernel.

The higher dimensional case was further studied Rochberg \cite{MR1358178}, and 
 Rochberg--Semmes \cites{MR845199,MR1021138}.  We will rely upon their innovations.  

\bigskip 

We will collect known results in \S \ref{s:2}, including Bloom's characterization and the Rochberg--Semmes notion of \emph{nearly weakly orthogonal.}  The next section 
\S \ref{s:BB} will address properties of the dyadic Besov spaces, and their relationwship with the continuous space.  
In the next two sections \S \ref{s:Schatten Pi} and \S \ref{s:Schatten Hilbert} 
we address the dyadic versions of these questions, where we can establish results for $1 < p < \infty$.  
Then, in \S \ref{s:Schatten Hilbert} we address the necessity for the Hilbert transform inequality.  Last of all in \S \ref{s: Besov further} we point to some future directions 
for research.

\begin{comment}
	\begin{remark}\label{coro schatten R}
We point out that it is natural to explore whether the result in Theorem  \ref{schatten H} holds for $S^p$ with $p\not=2$.  By \cite{MR500308}, a Russo type inequality, it seems to indicate that a suitable definition of $B_{\nu}^{p}(\mathbb R)$ could be the following:  ``Let $1<p<\infty$.  Suppose $\mu,\lambda\in A_2$ and set $\nu= \mu^{1\over2}\lambda^{-{1\over2}}$. Let $b\in L_{\rm loc}^{1}(\mathbb R)$. Then we say that $b$ belongs to the weight  Besov space $B_{\nu}^{p}(\mathbb R)$ if
\begin{align*}
{ \|b\|_{B_{\nu}^{p}(\mathbb R)}:= \Bigg( \int_{ \mathbb R} \int_{\mathbb R} {\big|b(x)-b(y)\big|^p\over |x-y|^2}\lambda^{p\over2}(x)\mu^{-{p\over2}}(y)dydx\Bigg)^{1\over p} }<\infty.''
\end{align*}

However, there are natural examples of $\lambda,\mu\in A_2$ such that $\lambda^{p\over2}$ and $\mu^{-{p\over2}}$ are not locally integrable. Thus, this $B_{\nu}^{p}(\mathbb R)$ does not connect to the dyadic version $B_{\nu,d}^{p}(\mathbb R)$ as in Definition \ref{def dy Bp}, and hence the approach via using dyadic structure does not work. More discussions will be provided in the last section.  The case $p\not=2$ is essentially different from $p=2$, which requires new idea and techniques.
\end{remark}

\end{comment}

%%%%%%%%%%%%%%%%%%%%%%%%%%%%%%%%%%
\section{Preliminaries}
\setcounter{equation}{0}
\label{s:2}

In this section we collect the necessary background for the theorems to be proved.
\subsection{Schatten Classes}
For two separable Hilbert spaces $\mathcal G$ and $\mathcal H$ over $\mathbb C$, we let
$\mathcal B(\mathcal G, \mathcal H)$ be the space of bounded linear operators from $\mathcal G$ to $\mathcal H$, and  $\mathcal B_0(\mathcal G, \mathcal H)$ the space of compact operators. Suppose
$T\in \mathcal B_0(\mathcal G, \mathcal H)$. Then $|T|: = (T^*T)^{1\over 2}$ is compact and positive.
The eigenvalues of $|T|$, which are non-negative, are called the singular values of $T$. Singular values are listed in decreasing order 
according to multiplicity by $ \sigma _1 \geq \sigma _2 \geq \cdots \geq0$.  For $0< p<\infty$,
the Schatten class $S^p(\mathcal G,\mathcal H)$ consists of those compact operators $T$ whose sequence of singular values (counted according to multiplicity and arranged in decreasing order) belong to $\ell^p$. In that case, we write 
\begin{equation} \label{e:SpNorm}
\lVert T \rVert_{S^p (\mathcal G, \mathcal H)} ^p = \sum _{j\geq 1} \sigma _j ^{p}. 
\end{equation}

In their innovative approach to singular values, Rochberg and Semmes \cite{MR1021138} introduced a notion of \emph{nearly weakly orthogonal} (NWO) sequences of functions. 
This is an important tool for us, so we recall facts here. We apply the notion in the unweighted section.  

\begin{definition} \label{d:NWO}  Let $1< p < \infty$, and  $ \mathcal D$ be a dyadic system, a collection of functions  $ \{ e_I \colon I\in \mathcal D\}$  is \emph{NWO}  if 
each $e_I$ is supported on $I$, and the maximal function below is bounded on $L^p ( \mathbb{R})$
\begin{equation}
\label{e:NWOdef}  \sup _I  \chi _I  \frac{ \lvert  \langle  f, e_I \rangle \rvert}{ \lvert  I \rvert ^{1/2}} . 
\end{equation}
\end{definition}
An elementary example of NWO sequence is any sequence of functions satisfying 
$ \lvert  e_ I\rvert \lesssim \lvert  I\rvert ^{-1/2} \chi _I $.  
Since we work with weights, the NWO sequences we meet will not satisfy this criteria.  
Instead, we will have need of this fact.

\begin{proposition}[ \cite{MR1021138}*{\S7D}] 
 \label{p:NWO} If the collection of functions $ \{e_I \colon  I\in \mathcal D\}$ are supported on $I$ and satisfy  for some $2< r < \infty $, 
$ \lVert e_I \rVert_{r} \lesssim \lvert  I \rvert ^{1/r-1/2}$, then $ \{e_I\} $ is NWO. 
\end{proposition}

The concept of NWO is related to Schatten norms by way of this  key inequality. For any 
bounded compact operator $ A $ on $ L ^2 (\mathbb R^{n})$:
\begin{equation}\label{e:NWO}
\sum_{I\in  \mathcal D } \lvert  \langle A e_I, f_I \rangle\rvert ^{p} \lesssim \lVert A \rVert _{S ^{p}} ^{p},
\end{equation}
where $\{e_I\}_{I\in \mathcal D}$ and $\{f_I\}_{ I \in \mathcal D}$ are  NWO. This inequality can be found in  \cite{MR1021138}*{(1.10), p. 242}.

\subsection{$A_2$ weights}

\begin{definition}%[$A_p$ weight]
  \label{def:Ap}
  Let $w(x)$ be a nonnegative locally integrable function
  on~$\R$. We
  say $w$ is an $A_2$ \emph{weight}, written $w\in
  A_2$, if
  \[
    [w]_{A_2}
    \coloneqq \sup_B {1\over|B|}\int_B w(x)dx
    \cdot {1\over |B|}\int_B
      w(x)^{-1}dx
    < \infty.
  \]
  Here the suprema are taken over all intervals $B\subset \R$.
  The quantity $[w]_{A_2}$ is called the \emph{$A_2$~constant
  of~$w$}.
\end{definition}

It is well known that $A_2$ weights are doubling. Namely,  
\begin{lemma}\label{lemcomparison}
Let $w\in A_2(\mathbb R)$. Then 
for every $\lambda>1$ and for every interval $B\subset \R$,
\begin{align}\label{doubling constant of weight}
w(\lambda B)  \lesssim \lambda \, w(B). 
\end{align}
See for example \cite{MR3243741}*{Lemma 9.2.1}.
\end{lemma}

\subsection{Weight BMO and VMO spaces}
We recall the  weighted BMO space and VMO space and their dyadic versions.
\begin{definition}\label{MWbmo}
Suppose $w \in A_\infty$.
A function $b\in L^1_{\rm loc}(\mathbb R)$ belongs to
${ BMO}_{w}(\mathbb R)$ if
\begin{equation*}
\|b\|_{{ BMO}_w(\mathbb R)}:=\sup_{B}{1\over w(B)}\int_{B}
|b(x)-b_{B}| \,dx<\infty,
\end{equation*}
where $b_{B}:= {1\over |B|}\int_B b(x)dx$ and the supremum is taken over all balls $B\subset \mathbb R$.  The dyadic weighted BMO space ${ BMO}_{w,d}(\mathbb R)$ consists of  functions $b\in L^1_{\rm loc}(\mathbb R)$
such that $\|b\|_{{ BMO}_{w,d}(\mathbb R)}<\infty$, where the $\|b\|_{{ BMO}_{w,d}(\mathbb R)}$ is defined the same as above with intervals replaced by dyadic cubes.
\end{definition}

\begin{definition}\label{vmo}
Suppose $w \in A_\infty$.
A function $b\in { BMO}_w(\mathbb R)$ belongs to
${ VMO}_w(\mathbb R)$ if
\begin{align*}
&{\rm (i)}\ \lim_{a\to0}\sup_{B:\  r_B=a}{1\over w(B)}\int_{B}
|b(x)-b_{B}| \,dx=0,\\
&{\rm (ii)}\ \lim_{a\to\infty}\sup_{B:\  r_B=a}{1\over w(B)}\int_{B}
|b(x)-b_{B}| \,dx=0,\\
&{\rm (iii)}\ \lim_{a\to\infty}\sup_{B\subset \mathbb R \backslash B(x_0,a)}{1\over w(B)}\int_{B}
|b(x)-b_{B}| \,dx=0,
\end{align*}
where $x_0$ is any fixed point in $\mathbb R$.  A  function $b\in { BMO}_{w,d}(\mathbb R)$ belongs to the dyadic weighted VMO space
${ VMO}_{w,d}(\mathbb R)$ if the above three limits hold with intervals replaced by dyadic cubes.

\end{definition}

\subsection{Operator Theory Lemma}
We provide a simple lemma that allows us to study  Schatten norms  on unweighted $L^2(\R)$.
\begin{lemma}\label{key lemma}
Suppose $1\leq p<\infty$ and $\mu,\lambda\in A_2$.
Then $T$ is in $S^p(L^2_\mu(\R) , L^2_\lambda(\R))$ if and only if
$\lambda^{1\over2}T\mu^{-{1\over2}}$ is in $S^p(L^2(\R) , L^2(\R))$. Moreover, we have
$$ \|T\|_{S^p(L^2_\mu(\R) , L^2_\lambda(\R))} = \|\lambda^{1\over2}T\mu^{-{1\over2}}\|_{S^p(L^2(\R) , L^2(\R))}. $$
\end{lemma}
\begin{proof}
We first point out that $T$ is bounded from $L^2_\mu(\R)$ to  $L^2_\lambda(\R)$ if and only if
$\lambda^{1\over2}T\mu^{-{1\over2}}$ is bounded from $L^2(\R)$ to $L^2(\R)$. Moreover, $T$ is compact from $L^2_\mu(\R)$ to  $L^2_\lambda(\R)$ if and only if $\lambda^{1\over2}T\mu^{-{1\over2}}$ is compact from $L^2(\R)$ to $L^2(\R)$.  The case of Schatten membership is equally straightforward.

To show the Schatten class estimate, we treat $\lambda^{1/2}$ as a multiplication operator, that is,
$$ M_{\lambda^{1\over2}}f(x) = \lambda^{1\over2}(x) f(x),\quad \forall f\in L^2_\lambda(\R).$$
Then we see that
$M_{\lambda^{1\over2}}$ is a unitary operator from $L_\lambda^2(\R)$ to $L^2(\R)$. Similarly, $M_{\mu^{-{1\over2}}}$ is a unitary operator from $L^2(\R)$ to $L_\mu^2(\R)$. Note that unitary operators preserve $S^p$ classes for $1\leq p<\infty$, we obtain our argument.
\end{proof}

\section{Bloom type Besov spaces and dyadic structure} \label{s:BB}

\subsection{Equivalent norms on the dyadic weighted Besov spaces}

Similarly to the weighted $BMO$ spaces, there are equivalent norms for the \emph{dyadic} Besov spaces defined in 
Definition \ref{def dy Bp}.   
\begin{proposition}\label{Bp properties}
Suppose $\mu,\lambda\in A_2$ and set $\nu= \mu^{1\over2}\lambda^{-{1\over2}}$. Then
\begin{align*}
\|b\|_{B_{\nu,d}^{p}(\mathbb R)} ^p
&\approx  \sum_{I\in\mathcal{D}} \Big|\widehat{b}(I) {\lambda(I)^{1\over2}  (\mu^{-{1}}(I))^{{1\over2}} \over |I|^{{3\over2}}}\Big|^p 
\\
&\approx  \sum_{I\in\mathcal{D}} \Big| \widehat{b}(I) {|I|^{1\over2}\over  {\lambda^{-1}(I)^{1\over2}  (\mu(I))^{{1\over2}}  }}\Big|^p 
\\
&\approx \sum_{I\in\mathcal D} \bigg( { |\widehat{b}(I)|\, |I|^{1\over2}\over {\nu(I)}} \bigg)^p,
\end{align*}
where the implicit constants depend only on $[\mu]_{A_2}$, $[\lambda]_{A_2}$ and $[\nu]_{A_2}$.
\end{proposition}

\begin{proof}
This depends upon elementary properties of $A_2$ weights.  In view of Definition \ref{def dy Bp}, it clearly suffices to argue that 
for any dyadic interval $I$, 
\begin{align} \label{e:AAprox}
\frac{\lvert  I \rvert } { \nu (I)}  & \approx 
{\lambda(I)^{1\over2}  (\mu^{-{1}}(I))^{{1\over2}} \over |I|} 
\approx 
{|I|\over  {\lambda^{-1}(I)^{1\over2}  (\mu(I))^{{1\over2}}  }} . 
\end{align}
But, the definition of the $\sigma \in A_2$ is  
$$ { \sigma(I)\over |I| }{ \sigma^{-1}(I)\over |I| } \approx 1 .$$
It follows that 
\begin{equation}
{\lambda(I) \mu^{-{1}}(I) \over |I| \cdot \lvert  I \rvert }
 \approx 
{|I| \cdot \lvert  I \rvert \over  {\lambda^{-1}(I)  \mu(I)  }}. 
\end{equation}
This is the rightmost equivalence in \eqref{e:AAprox}. 

Recalling the definition of $\nu $ and using Cauchy--Schwarz,  we see that
\begin{align*}
 \frac{ \lvert  I \rvert ^2  }{ \nu (I)} 
 & \approx \nu^{-1}(I)  
\\
 &=\int_I \nu^{-1}(x)dx = \int_I \mu^{-{1\over2}}(x)\lambda^{{1\over2}}(x) dx
 \\
 & \leq \mu^{-1}(I)^{1\over2} \lambda(I)^{1\over2}.
\end{align*}
And, again by Cauchy-Schwarz, 
\begin{align*}
{|I|\over  {\lambda^{-1}(I)^{1\over2}  (\mu(I))^{{1\over2}}  }}  
& \lesssim  \frac{ \lvert I \rvert}{ v (I)},
\end{align*}
which completes the proof. 
\end{proof}

\subsection{Dyadic structure of the weighted Besov space $B_{\nu}^{2}(\mathbb R)$ }

In this subsection we will prove Theorem \ref{thm dyadic structure}, namely that the continuous Besov space is the intersection of well chosen dyadic Besov spaces.  
There are two halves. In the first half, the  continuous Besov space $B_{\nu}^{2}(\mathbb R)$ is contained a  dyadic Besov space $B_{\nu,d}^{2}(\mathbb R)$, 
for any choice of dyadic system $ \mathcal D$.   To be more explicit, we have the following.

\begin{proposition}\label{prop B and Bd}
Suppose $\mu,\lambda\in A_2$ and set $\nu= \mu^{1\over2}\lambda^{-{1\over2}}$.
Then we have $B_{\nu}^{2}(\mathbb R)\subset B_{\nu,d}^{2}(\mathbb R)$, where $B_{\nu,d}^{2}(\mathbb R)$ is the dyadic weighted Besov space associated with an arbitrary dyadic system ${\mathcal D}\subset \mathbb R$ and we have
$$ \|b\|_{B_{\nu,d}^{2}(\mathbb R)}\lesssim \|b\|_{B_{\nu}^{2}(\mathbb R)}$$
with the implicit constants depending only on $[\mu]_{A_2}$ and $[\lambda]_{A_2}$.
\end{proposition}
\begin{proof}
It suffices to see that the dyadic Besov norm is dominated by the continuous Besov norm.   The key inequalities here are specific to a choice of dyadic interval $I$. 
Below, let $I' = I + 2\lvert  I \rvert$.  We will freely use the $A_2$ property. 
\begin{align*}
{ \big|\int_Ib(x)h_I(x)dx\big|\, |I|^{1\over2}\over {\nu(I)}} 
& \lesssim \bigg|\int_Ib(x)h_I(x)dx\bigg|\, {\nu^{-1}(I)\over |I|^{3\over2}}
\\
& \lesssim\inf_{y\in   I'}  \bigg|\int_I\Big(b(x)-b(y)\Big)h_I(x)dx\bigg|\, {\nu^{-1}( I')\over |I|^{3\over2}}
\\
& \lesssim\int_{  I'}  \bigg|\int_I\Big(b(x)-b(y)\Big)h_I(x)dx\bigg|\, \frac{\nu^{-1}(y)} {|I|^{3\over2}} dy 
\\
& \lesssim \int_{  I'}  \int_I\Big|b(x)-b(y)\Big||h_I(x)| \lambda^{1\over2}(x)\lambda^{-{1\over2}}(x)dx\,\lambda^{1\over2}(y)\mu^{-{1\over2}}(y)dy \ {1\over |I|^{3\over2}} 
\\
& \lesssim  \bigg(\int_{  I'} \int_I\big|b(x)-b(y)\big|^2\lambda(x)\mu^{-1}(y)dydx  \times \int_{ I'} \int_I |h_I(x)|^2\lambda^{-1}(x)\lambda(y)dydx\bigg)^{1\over2}\, \ {1\over |I|^{3\over2}} 
\\
& \lesssim \bigg(\int_{ I'} \int_I{\big|b(x)-b(y)\big|^2\over |x-y|^2}\lambda(x)\mu^{-1}(y)dydx \ {\lambda(I)\lambda^{-1}(I) \over |I|^2} \bigg) ^{\frac12} 
 \\
 & \lesssim \bigg(\int_{ I'} \int_I{\big|b(x)-b(y)\big|^2\over |x-y|^2}\lambda(x)\mu^{-1}(y)dydx \bigg) ^{\frac12} . 
\end{align*}
Square this, and sum over $I\in \mathcal D$, to get an expression dominated by $\lVert b \rVert_{B ^2 _ \nu (\R)} ^2$.  
The proof of Proposition \ref{prop B and Bd} is complete.
\end{proof}

The second half of the proof is to show that the continuous Besov space is in the intersection of two well chosen dyadic Besov spaces.  
The dependence of the dyadic Besov space on the choice of dyadic system is supressed in the notation. For the purposes 
of this next proposition,  we use the notation  $B ^2_ \nu ( \mathbb{R}, \mathcal D)$ for the dyadic Besov space associated with dyadic system $\mathcal D$.  

\begin{proposition}\label{prop B and Bd 2}
Suppose $\mu,\lambda\in A_2$ and set $\nu= \mu^{1\over2}\lambda^{-{1\over2}}$. 
There are two choices of grids $ \mathcal D^0$ and $ \mathcal D^1$ for which we have  
 $ B_{\nu}^{2}(\mathbb R , \mathcal D^0)\cap  B_{\nu}^{2}(\mathbb R, \mathcal D^1) = B_{\nu}^{2}(\mathbb R)$. 
Moreover, we have 
 $$\|b\|_{B_{\nu}^{2}(\mathbb R)}\approx   \|b\|_{B_{\nu}^{2}(\mathbb R , \mathcal D^0)}+ \|b\|_{B_{\nu}^{2}(\mathbb R , \mathcal D^1)}.$$
\end{proposition}

\begin{proof}
 We take two dyadic grids $ \mathcal D^0 $ and $\mathcal D^1$ so that for all intervals $I$ there is a $Q\in 
 \mathcal D^0 \cup \mathcal D^1$ with 
 \begin{equation}
 \label{e:oneThird}   I\subset Q \subset 4 I . 
 \end{equation}
One option is that  $ \mathcal D^0 $ is the standard dyadic system in $\R$ and  $ \mathcal D^1 $ is the  `'one-third shift' 
 of  $ \mathcal D^0 $, see for example \cite{MR3420475}.

By Proposition \ref{prop B and Bd}, it  suffices to show that the continuous Besov norm is dominated by the larger of the two dyadic Besov space norms.  Below, $\mathcal D$ is the standard dyadic grid in $\R$. 
Also,  $a , n \in \mathbb N $ are integers that we will specify.  Observe that:
\begin{align}
 \|b\|_{B_{\nu}^{2}(\mathbb R)} ^2
& =   \int_{ \mathbb R} \int_{\mathbb R} {\big|b(x)-b(y)\big|^2\over |x-y|^2}\lambda(x)\mu^{-1}(y)dydx
\\
&=
 \sum_{I\in\mathcal D}\int_{I} \int_{\{y\in \mathbb R:  2 ^{a} \lvert  I \rvert<|x-y|\leq 2 ^{a+1} \lvert  I \rvert}{\big|b(x)-b(y)\big|^2\over |x-y|^2}\lambda(x)\mu^{-1}(y)dydx
\\
&=
 \sum_{I\in\mathcal D} \lvert  I \rvert ^{-2}
 \int_{I} \int_{\{y\in \mathbb R:  2 ^{a} \lvert  I \rvert<|x-y|\leq 2 ^{a+1} \lvert  I \rvert}{\big|b(x)-b(y)\big|^2 }\lambda(x)\mu^{-1}(y)dydx
\\  \label{e:IA}
& \lesssim 
 \sum_{I\in\mathcal D} \lvert  I \rvert ^{-2} \sum _{m=n} ^{2n-1}
 \int_{I} \int_{\{y\in \mathbb R:  2 ^{a} (m/n)\lvert  I \rvert<|x-y|\leq 2 ^{a} ((m+1)/n) \lvert  I \rvert} {\big|b(x)-b(y)\big|^2 }\lambda(x)\mu^{-1}(y)dydx
\end{align}
The second integral above is over a symmetric interval. Consider the two intervals 
\begin{equation} \label{e:IAK}
I, \qquad    I + 2 ^{a}\lvert  I \rvert[m/n, (m+1)/n],  \qquad  n\leq m < 2n.  
\end{equation}
Now, we  choose $a=5$, and $n =1000$, so the second interval is smaller in length, but still comparable to $I$ in length.  And, they are separated by a distance approximately 
$2^a \lvert  I \rvert$. By \eqref{e:oneThird},  we can choose $a$ so that there is a dyadic $I'\in \mathcal D ^{0} \cup \mathcal D ^{1}$ 
which contains both intervals above, and moreover $I$ is contained in the left half of $I'$, and $I + 2 ^{a}\lvert  I \rvert[m/n, (m+1)/n]$ the right half. 
We can argue similarly for $I -2 ^{a}\lvert  I \rvert[m/n, (m+1)/n]$. Below we continue with $I$ and $I'$. 
In particular, for the main term in  \eqref{e:IA}, with fixed $n\leq m < 2n$,  we have 
\begin{align}
\lvert  I \rvert ^{-2} 
 \int_{I} &
 \int_{\{y \colon   2 ^{a} (m/n)\lvert  I \rvert<|x-y|\leq 2 ^{a} ((m+1)/n) \lvert  I \rvert}
  {\big|b(x)-b(y)\big|^2 }\lambda(x)\mu^{-1}(y)dydx
\\&= 
\lvert  I \rvert ^{-2}\int_{I} 
\int_{\{y \colon   2 ^{a} (m/n)\lvert  I \rvert<|x-y|\leq 2 ^{a} ((m+1)/n) \lvert  I \rvert}
{\big|b(x)-B_{I'} +B_{I'} -b(y)\big|^2}\lambda(x)\mu^{-1}(y)dydx 
\\
& \lesssim  
 \lvert  I \rvert^{-2} \Biggl\{ \int_{I'} \int_{I'}{\big|b(x)-B_{I'}\big|^2}\lambda(x)\mu^{-1}(y)dydx 
 +
 \int_{I'} \int_{I'}{\big|b(y)-B_{I'}\big|^2}\lambda(x)\mu^{-1}(y)dydx 
 \Biggr\}.
\end{align}
Above, we choose $B_{I'} = \lvert  I' \rvert ^{-1}\int _{I'} b (t)\;dt$.  

It follows that the norm  $ \lVert b \rVert_{B ^2_ \nu (\mathbb{R})} ^2$ is dominated by several terms, one of which is 
\begin{equation} \label{e:1of8}
 \sum_{I\in\mathcal D^0} 
 \lvert  I\rvert^{-2} \int_{I} \int_{I}{\big|b(x)-B_I\big|^2}\lambda(x)\mu^{-1}(y)dydx. 
\end{equation}
The other terms are obtained by varying the role of $m$ in \eqref{e:IAK}, considering the negative of the  intervals in \eqref{e:IAK}, 
exchanging the role of the dyadic grid, and  the roles of $x$ and $y$. All cases are similar, so we continue with the one above.  
In \eqref{e:1of8}, the point is that 
\begin{equation}
\mathbf 1 _{I} (b(x)-B_I) = \sum _{\substack{J\in \mathcal D ^{0} \\  J\subset I}} \langle  b,h_J \rangle h_J. 
\end{equation}
That is, only the smaller scales contribute.  But then, it is straight forward to see that we can make a pure sum on scales.  
\begin{align}
\eqref{e:1of8} & \lesssim 
\sum_{I\in\mathcal D^0} 
 \lvert  I\rvert^{-2} \int_{I} \int_{I}{\big|\langle  b,h_I \rangle h_I\big|^2}\lambda(x)\mu^{-1}(y)dydx 
 \lesssim  \lVert b \rVert_{B ^2 _{\nu } (\mathbb{R}, \mathcal D ^{0})} ^2. 
\end{align}
The last inequality appeals to Proposition \ref{Bp properties}. This completes the proof.  
\end{proof}

The proof of Theorem \ref{thm dyadic structure} is complete.

\subsection{Comparison of $B_{\nu}^{2}(\mathbb R)$ with weighted BMO and VMO spaces}

We now recall the dyadic weighted BMO space ${BMO}_{\nu,d}(\mathbb R)$ and give equivalent norms for it.

\begin{lemma}\label{Bp and VMOd}
Suppose that $\nu\in A_2$. Then we have
$$B_{\nu,d}^{2}(\mathbb R)\subset {\rm VMO}_{\nu,d}(\mathbb R).$$
\end{lemma}
\begin{proof}
We first prove that $B_{\nu}^{2}(\mathbb R)\subset {\rm BMO}_{\nu,d}(\mathbb R).$ 
Recall that 
\begin{align*}
\|b\|_{{ BMO}_{\nu,d}(\mathbb R)}^2\approx \sup_{K\in\mathcal D}{1\over \mu^{-1}(K)}\sum_{ I\subset K} |\widehat{b}(I)|^2 {\mu^{-1}(I)^2 \over |I|^2 } {\lambda(I)\over |I|}.  \end{align*}

For any $K\in\mathcal D$, by Proposition \ref{Bp properties}, we have
\begin{align*}
{1\over \mu^{-1}(K)}\sum_{ I\subset K} |\widehat{b}(I)|^2 {\mu^{-1}(I)^2 \over |I|^2 } {\lambda(I)\over |I|}
&\leq \sum_{ I\subset K} |\widehat{b}(I)|^2 {\mu^{-1}(I) \over |I|^2 } {\lambda(I)\over |I|}  \lesssim \|b\|_{B_{\nu}^{2}(\mathbb R)}^2.
\end{align*}
Taking the supremum on the left-hand side over $K$, we get that $ \|b\|_{{\rm BMO}_{\nu,d}(\mathbb R)}^2\lesssim \|b\|_{B_{\nu}^{2}(\mathbb R)}^2 $.  Hence, we see that $B_{\nu}^{p}(\mathbb R)\subset {\rm BMO}_{\nu,d}(\mathbb R)$. Moreover, since
$$  \|b\|_{B_{\nu}^{2}(\mathbb R)}^2\approx \sum_{ I\in\mathcal D} {1\over \mu(I)  } |\widehat{b}(I)|^2  {\lambda(I)\over |I|} <\infty,$$
it is direct to see that
\begin{align*}
&{\rm (i)}\ \lim_{a\to0} \sum_{ I\in\mathcal D, |I|<a} {1\over \mu(I)  } |\widehat{b}(I)|^2  {\lambda(I)\over |I|} =0,\\
&{\rm (ii)}\ \lim_{a\to\infty} \sum_{ I\in\mathcal D, |I|>a} {1\over \mu(I)  } |\widehat{b}(I)|^2  {\lambda(I)\over |I|} =0,\\
&{\rm (iii)}\ \lim_{a\to\infty} \sum_{ I\in\mathcal D, I \subset \R\backslash B(x_0,a)} {1\over \mu(I)  } |\widehat{b}(I)|^2  {\lambda(I)\over |I|} =0,
\end{align*}
where $x_0$ is any fixed point in $\mathbb R$.  Hence, $b\in  {\rm VMO}_{\nu,d}(\mathbb R)$.
\end{proof}

As a corollary of Lemma \ref{Bp and VMOd}, Proposition \ref{prop B and Bd} and the relationship of weighted VMO and weighted dyadic VMO spaces via intersection \cite{MR3420475}
 we see that for every $\nu\in A_2$,
$$B_{\nu}^{2}(\mathbb R)\subset {\rm VMO}_{\nu}(\mathbb R).$$
This, together with the two weight compactness argument of $[b,H]$ in \cite{MR4345997}, gives that:
for $\mu,\lambda\in A_2$ and  $\nu= \mu^{1\over2}\lambda^{-{1\over2}}$,
for every $b\in B_{\nu}^{2}(\mathbb R)$, $[b,H]$ is compact from $L^2_\mu(\R)$ to $L^2_\lambda(\R)$.

\bigskip
\bigskip

\section{Schatten classes and paraproducts}
\label{s:Schatten Pi}
\setcounter{equation}{0}

 Let ${\mathcal D}$ be an arbitrary dyadic system in $ \mathbb R$.
Recall that  $\Pi_b\equiv \sum_{I\in\mathcal{D}} \widehat{b}(I) h_I\otimes \frac{\mathsf{1}_I}{\left\vert I\right\vert}$ is the `paraproduct' operator with symbol function $b$. The adjoint of $\Pi_b$ on \textit{unweighted} $L^2(\mathbb{R})$ is defined by
$\Pi_b^{\ast} \equiv \sum_{I\in\mathcal{D}} \widehat{b}(I) \frac{\mathsf{1}_I}{\left\vert I\right\vert}\otimes h_I$,
where $\widehat{b}(I) =\langle b,h_I\rangle$. Moreover, we see that for $\mu,\lambda\in A_2$,
\begin{equation}
\lVert \Pi _b   \rVert_{L^2_\mu(\R) \rightarrow L^2_\lambda(\R)} 
=
\lVert \Pi _b  ^{\ast} \rVert_{L^2_{ \lambda ^{-1}}(\R)  \rightarrow L^2_{ \mu ^{-1}} (\R)}.
\end{equation}

We have the following characterization of the Schatten membership of the paraproduct operators.

\begin{lemma}\label{Pi Sp norm}
Suppose $1<p<\infty$, $\mu,\lambda\in A_2$ and set $\nu= \mu^{1\over2}\lambda^{-{1\over2}}$.
 Suppose   $b\in {\rm VMO}_{\nu,d} (\mathbb R)$. Then we have $\Pi_b$ is in $S^p(L^2_\mu(\R) , L^2_\lambda(\R))$ if and only if $b$ is in $B_{\nu,d}^p(\R)$.  Moreover,
$$\|\Pi_b\|_{ S^p(L^2_\mu(\R) , L^2_\lambda(\R)) } \approx \|b\|_{B_{\nu,d}^p(\R)}.$$
\end{lemma}

\begin{proof}
\noindent
{\bf Sufficiency:}

Suppose $b\in B_{\nu, d}^p(\R)$.  By definition we see that
\begin{align*}
(\lambda^{1\over 2} \Pi_b \mu^{-{1\over2}})(f)(x) &=\sum_{I\in\mathcal{D}} \widehat{b}(I) \lambda^{1\over 2}(x)h_I(x) \int_{\mathbb{R}} \frac{\mathsf{1}_I(y)}{\left\vert I\right\vert} f(y)\,  \mu^{-{1\over2}}(y)\, dy\\
&=\sum_{I\in\mathcal{D}} \widehat{b}(I) {\lambda(I)^{1\over2}  (\mu^{-{1}}(I))^{{1\over2}} \over |I|^{1\over2}} {1\over |I|} \cdot  {\lambda^{1\over 2}(x) |I|^{1\over2}h_I(x)\over \lambda(I)^{1\over2} }  \int_{\mathbb{R}}  f(y)\,   \frac{\mathsf{1}_I(y) \mu^{-{1\over2}}(y)}{ (\mu^{-{1}}(I))^{{1\over2}} } \, dy\\
&= \sum_{I\in\mathcal{D}} B(I) \cdot  G_I(x) \int_{\mathbb{R}}  f(y)\,  H_I(y) dy,
\end{align*}
where
$$ B(I)= \widehat{b}(I) {\lambda(I)^{1\over2}  (\mu^{-{1}}(I))^{{1\over2}} \over |I|^{{3\over2}}}, \quad G_I(x)={\lambda^{1\over 2}(x)|I|^{1\over2} h_I(x)\over \lambda(I)^{1\over2}  },\quad H_I(y)= \frac{\mathsf{1}_I(y) \mu^{-{1\over2}}(y)}{ (\mu^{-{1}}(I))^{{1\over2}} } .$$
The purpose of these equalities is that the last line places us in the position of using the NWO inequality \eqref{e:NWO} to conclude the proof.   
This depends upon us showing that $\{ G_I\}$ and $ \{H_I\}$ are NWO for $L^2$, in the language of Definition \ref{d:NWO}.  
It then follows from \eqref{e:NWO} that 
\begin{align*}
  \|\Pi_b\|_{ S^p(L^2_\mu(\R) , L^2_\lambda(\R)) } ^{p}  &= \|\lambda^{1\over 2} \Pi_b \mu^{-{1\over2}}\|_{ S^p( L^2(\R) , L^2(\R) ) } ^{ p} 
  \\
& 
 \lesssim \sum_{I\in\mathcal{D}} | B(I)|^p \approx  \|b\|_{ B_{\nu}^p(\R)} ^{p}.
\end{align*}

The arguments that $\{ G_I\}$ and $ \{H_I\}$ are NWO depend only on the reverse H\"older inequality for $A_2$ weights. They are similar, and we address $G_I$, 
using the sufficient condition from Proposition \ref{p:NWO}, which requires us to control a local norm in an index $r >2$.  
To see this, we first note that $G_I$ is supported in $I$. Next, as  $ \lambda $ is $A_2$, there is a   $r = r _{[ \lambda ]_ {A_2}}>2$, so 
that for each dyadic interval $I$ 
\begin{equation*}
\Bigl[  \lvert  I \rvert ^{-1}\int _I  \lambda ^{r/2}  \; dx \Bigr] ^{2/r} \lesssim \frac{ \lambda (I)}{ \lvert  I \rvert}. 
\end{equation*}
With this choice of $r$, we can estimate as follows.
\begin{align*}
\|G_I\|_{L^r(\R^d)}&\leq  {|I|^{1\over2}\over \lambda(I)^{1\over2}}\bigg( \int_{I} {\lambda^{r\over 2}(x)|h_I(x)|^r  } dx\bigg)^{1\over r}\leq {|I|^{1\over2}\over \lambda(I)^{1\over2}}\bigg( {|I|\over |I|^{r\over2}}\ {1\over |I|} \int_{I} {\lambda^{r\over 2}(x)  } dx\bigg)^{1\over r}\\
&= {|I|^{1\over2}\over \lambda(I)^{1\over2}}\bigg( {|I|\over |I|^{r\over2}}\bigg)^{1\over r}\bigg( {1\over |I|} \int_{I} {\lambda^{r\over 2}(x)  } dx\bigg)^{{2\over r}\cdot {1\over2} }\\
&\lesssim {|I|^{1\over2}\over \lambda(I)^{1\over2}}\bigg( {|I|\over |I|^{r\over2}}\bigg)^{1\over r}\bigg( {1\over |I|} \int_{I} {\lambda(x)  } dx\bigg)^{ {1\over2} }\\
&= |I|^{{1\over r}-{1\over 2}},
\end{align*}
where the last inequality follows from reverse H\"older inequality for $\lambda\in A_2$. Hence, the claim holds.

\medskip
\noindent {\bf Necessity:}
Suppose $\Pi_b \in S^p( L^2_\mu(\R) , L^2_\lambda(\R) )$. Then from Lemma \ref{key lemma} we have that $\lambda^{1\over 2} \Pi_b \mu^{-{1\over2}} \in S^p( L^2(\R) , L^2(\R) )$.  Note that for each dyadic interval $I$, we have $\widehat{b}(I)  = \langle \Pi_b(\chi_I),h_I \rangle$.  Thus,
\begin{align*}
\sum_{I\in\mathcal{D}} \Big|\widehat{b}(I) {\lambda(I)^{1\over2}  (\mu^{-{1}}(I))^{{1\over2}} \over |I|^{{3\over2}}}\Big|^p &=\sum_{I\in\mathcal{D}} \Big|\langle \Pi_b(\chi_I),h_I \rangle  { |I|^{1\over2}\over \lambda^{-1}(I)^{1\over2}  (\mu(I))^{{1\over2}}  }\Big|^p\\
&=\sum_{I\in\mathcal{D}} \Big|\big\langle \lambda^{1\over 2} \Pi_b \mu^{-{1\over2}}(\mu^{{1\over2}}\chi_I),\lambda^{-{1\over 2}}h_I \big\rangle  { |I|^{1\over2}\over \lambda^{-1}(I)^{1\over2}  (\mu(I))^{{1\over2}}  }\Big|^p\\
&=\sum_{I\in\mathcal{D}} \big|\big\langle \lambda^{1\over 2} \Pi_b \mu^{-{1\over2}}G_I
, H_I\big\rangle  \big|^p.
\end{align*}
Above, we use the notation 
$$
G_I:={\mu^{{1\over2}}\chi_I\over   \mu(I)^{{1\over2}}} \quad{\rm and}\quad H_I:={ \lambda^{-{1\over 2}}h_I |I|^{1\over2} \over \lambda^{-1}(I)^{1\over2}  }.
$$
These two collections of functions are NWO, in the sense of Definition \ref{d:NWO}. 
(We argued along this lines in the sufficiency direction.)  
  And so we have by \eqref{e:NWO} that
\begin{align*}
\sum_{I\in\mathcal{D}} \Big|\widehat{b}(I) {\lambda(I)^{1\over2}  (\mu^{-{1}}(I))^{{1\over2}} \over |I|^{{3\over2}}}\Big|^p
\lesssim \|\lambda^{1\over 2} \Pi_b \mu^{-{1\over2}}\|_{ S^p( L^2(\R) , L^2(\R) ) } ^{p}. 
\end{align*}
That is, we have established that 
$ \|b\|_{ B_{\nu, d}^p(\R)}\lesssim \|\lambda^{1\over 2} \Pi_b \mu^{-{1\over2}}\|_{ S^p( L^2(\R) , L^2(\R) ) }$. 
\end{proof}

Based on Lemma \ref{Pi Sp norm}  we have the following.
\begin{coro}\label{Pi* Sp norm}
Suppose $1<p<\infty$, $\mu,\lambda\in A_2$ and set $\nu= \mu^{1\over2}\lambda^{-{1\over2}}$.
 Suppose   $b\in {\rm VMO}_{\nu} (\mathbb R)$. Then we have $\Pi_b\in S^p(L^2(\lambda^{-1}), L^2(\mu^{-1}))$ if and only if $b \in B_{\nu,d}^p(\R^d)$.  Moreover,
$$\|\Pi_b\|_{ S^p(L^2(\lambda^{-1}) , L^2(\mu^{-1})) } \approx \|b\|_{B_{\nu,d}^p(\R)}.$$
Hence, we have
$$\|\Pi_b^*\|_{ S^p(L^2_\mu(\R) , L^2_\lambda(\R)) } \approx \|b\|_{B_{\nu,d}^p(\R)}.$$
\end{coro}

\medskip
\section{Schatten class and commutator of Haar multiplier}
\setcounter{equation}{0}
\label{s:Schatten Haar}

 Let ${\mathcal D}$ be an arbitrary dyadic system in $ \mathbb R$.
Recall that the Haar multiplier is defined  in terms of an arbitrary function $ \epsilon \colon \mathcal D \to \{ -1, 1\} $ with
\begin{equation*}
T_\epsilon=  \sum _{I\in \mathcal D} \epsilon_I \widehat f (I) h_I. 
\end{equation*}

We have the following characterization of the Schatten membership of the paraproduct operators.

\begin{lemma}\label{Para Sp norm}
Suppose $1<p<\infty$, $\mu,\lambda\in A_2$ and set $\nu= \mu^{1\over2}\lambda^{-{1\over2}}$.
 Suppose   $b\in {\rm VMO}_{\nu,d} (\mathbb R)$. Then we have $[b,T_\epsilon]$ is in $S^p(L^2_\mu(\R) , L^2_\lambda(\R))$ if and only if $b$ is in $B_{\nu,d}^p(\R)$.  Moreover,
$$\|[b,T_\epsilon]\|_{ S^p(L^2_\mu(\R) , L^2_\lambda(\R)) } \approx \|b\|_{B_{\nu,d}^p(\R)}.$$
\end{lemma}

\begin{proof}

{\bf Sufficiency:} Suppose $b\in B_{\nu,d}^p(\R)$.

Consider the commutator of Haar multiplier $[b,T_\varepsilon]$. We have
\begin{align*}
[b,T_\varepsilon] = \Pi_b T_\varepsilon-T_\varepsilon\Pi_b + \Pi^*_bT_\varepsilon-T_\varepsilon \Pi_b^*.
\end{align*}

Then it is direct that
\begin{align*}
\|[b,T_\varepsilon] \|_{ S^p(L^2_\mu(\R) , L^2_\lambda(\R)) } &\leq \|\Pi_b\|_{ S^p(L^2_\mu(\R) , L^2_\lambda(\R)) }\| T_\varepsilon\|_{L^2_\mu(\R)\to L^2_\mu(\R)}+\|T_\varepsilon\|_{L^2_\lambda(\R)\to L^2_\lambda(\R)}\|\Pi_b\|_{ S^p(L^2_\mu(\R) , L^2_\lambda(\R)) } \\
&\quad+ \|\Pi^*_b\|_{ S^p(L^2_\mu(\R) , L^2_\lambda(\R)) }\|T_\varepsilon\|_{L^2_\mu(\R)\to L^2_\mu(\R)}+\|T_\varepsilon\|_{L^2_\lambda(\R)\to L^2_\lambda(\R)} \|\Pi_b^*\|_{ S^p(L^2_\mu(\R) , L^2_\lambda(\R)) }.
\end{align*}

\bigskip

Based on the result in Section \ref{s:Schatten Pi}, 
we see that
$[b,T_\varepsilon]\in { S^p(L^2_\mu(\R) , L^2_\lambda(\R)) }$.
\bigskip
\bigskip

{\bf Necessity:} Suppose that
$[b,T_\varepsilon]\in { S^p(L^2_\mu(\R) , L^2_\lambda(\R)) }$.  Then we know that $ \lambda^{1\over2}[b,T_\varepsilon]\mu^{-{1\over2}}\in { S^p(L^2(\R) , L^2(\R)) }$.  Based on the non-degenerated condition on the Haar multiplier coefficients $\{\epsilon_I\}_I$, by using the same proof as that for Proposition \ref{prop schattenlarge1}, we see that the necessity holds.
\end{proof}

\bigskip

%%%%%%%%%%%%%%%%%%%%%%%%%%%%%%%%%%%%%%
\section{Schatten classes and commutator of the Hilbert transform}
\setcounter{equation}{0}
\label{s:Schatten Hilbert}

In this section we prove Theorem \ref{schatten H}. To show the sufficiency, we use the dyadic result together with the observation of Petermichl \cite{MR1756958}, that the Hilbert transform  is an average of Haar shifts.  For the necessity, we use the pointwise kernel lower bound and the median to split the Besov norm into the product of the commutator against NWOs.

We now provide the proof of our main theorem.

\smallskip
\noindent
{\bf Sufficiency:}
\begin{proposition}\label{prop schattenlarge0}
Suppose $\mu,\lambda\in A_2$ and set $\nu= \mu^{1\over2}\lambda^{-{1\over2}}$.
Suppose that $b\in B_{\nu}^{2}(\R)$, then
\begin{align*}
\|[b,H]\|_{S^2(L^2_\mu(\R) , L^2_\lambda(\R)) } \lesssim \|b\|_{B_{\nu}^{2}(\R)}.
\end{align*}
\end{proposition}
\begin{proof}

Suppose that $b\in B_{\nu}^{2}(\R)$. Then we know that $b\in B^2_{\nu,d}(\R)$ for any system of dyadic intervals. That is
$$ \sum_{I\in\mathcal{D}} \Big|\widehat{b}(I) {\lambda(I)^{1\over2}  (\mu^{-{1}}(I))^{{1\over2}} \over |I|^{{3\over2}}}\Big|^p <\infty.$$

We use Petermichl's observation that the Hilbert transform can be recovered through an appropriate average of Haar shifts, \cite{MR1756958}. On  a fixed dyadic lattice $\mathcal{D}$ with Haar basis $\{h_I\}_{I\in\mathcal{D}}$, we define a Haar shift operator  $\Sh$, following Petermichl.
 Set for any Haar function $h_I$,  
$\Sh h_I=\frac{1}{\sqrt{2}}(h_{I_-}-h_{I_+})$, which uniquely defines the operator. 
  Then, the Hilbert transform is an average of shift operators, with the average performed over the class of
all dyadic grids.  In particular, to prove norm inequalities for the Hilbert transform, it suffices to prove them
for the Haar shift operator.

The commutator with the Haar shift operator has an explicit expansion in terms of the paraproducts and $\Sh $,  
 \begin{equation}
\label{e:expand}
[b,\Sh]f=\Sh(\Pi_bf)-\Pi_b(\Sh f)+\Sh (\Pi_b^* f)-\Pi_b^*(\Sh f)+\Pi_{\Sh f} b-\Sh(\Pi_f b).
\end{equation}
It is known that for any $w\in A_2$, $\left\Vert \Sh\right\Vert_{L^2(w)\to L^2(w)}\lesssim 1$.  
There are six terms on the right above.  The first four are of the form a paraproduct (or its adjoint) composed with a bounded operator.  
But, the Schatten classes are closed under composition with bounded operators.  That is, we have, for instance for the first term on the right in \eqref{e:expand}, 
\begin{equation*}
\lVert \Sh(\Pi_bf) \rVert_{ S^2(L^2_\mu(\R) , L^2_\lambda(\R)) } \leq  \|\Sh\|_{L^2(\lambda)\to L^2(\lambda)} \|\Pi_b\|_{ S^2(L^2_\mu(\R) , L^2_\lambda(\R)) }
\lesssim  \|b\|_{B^2_{\nu}(\R)}. 
\end{equation*}
We have to appeal to Lemma \ref{Pi Sp norm} and Corollary  \ref{Pi* Sp norm}. The discussion of the next three terms is the same.  

The last two terms on the right in \eqref{e:expand} are not of this form, and are treated separately.
Set 
\begin{align*}
 \mathcal R_b f & \coloneqq \Pi_{\Sh f} b-\Sh(\Pi_f b) 
 \\
& =\sum_{I\in\mathcal{D}} \frac{\widehat{b}(I)}{\left\vert I\right\vert^{\frac{1}{2}}} \widehat{f}(I) (h_{I_+}(x)-h_{I_-}(x))
\\
& \eqqcolon \sum_{I\in\mathcal{D}} \frac{\widehat{b}(I)}{\left\vert I\right\vert^{\frac{1}{2}}} \widehat{f}(I) k_I(x).
\end{align*}
Expressed in this way, we see that the remainder is better behaved than the paraproducts, since the $k_I$ have strong orthogonality properties.  A straight forward variant of the proof of the Schatten norm estimates for paraproducts will complete the proof. The details are below.   

Recall from Lemma \ref{key lemma} that $\mathcal R_b \in S^2(L^2_\mu(\R) , L^2_\lambda(\R))$ if and only if $ \lambda^{1\over 2} \mathcal R_b \mu^{-{1\over2}}\in S^2( L^2(\R) , L^2(\R) )$.  We have that
\begin{align*}
(\lambda^{1\over 2} \mathcal R_b \mu^{-{1\over2}})(f)(x)&=  \sum_{I\in\mathcal{D}} \frac{\widehat{b}(I)}{\left\vert I\right\vert^{\frac{1}{2}}}  \int_{\mathbb R} f(y) \mu^{-{1\over2}}(y) h_I(y)dy  \lambda^{1\over 2}(x)k_{I}(x)
\\
&=\sum_{I\in\mathcal{D}} \widehat{b}(I) {\lambda(I)^{1\over2}  (\mu^{-{1}}(I))^{{1\over2}} \over |I|^{{3\over2}}}  \int_{\mathbb R} f(y) {\mu^{-{1\over2}}(y) |I|^{1\over2} h_I(y)\over  (\mu^{-{1}}(I))^{{1\over2}}}\,dy \ {\lambda^{1\over 2}(x) |I_+|^{1\over2}k_{I}(x)\over \lambda(I)^{1\over2}}\\
&=\sum_{I\in\mathcal{D}} B(I) \int_{\mathbb R} f(y) H_I(y)\,dy \ K_I(x),
\end{align*}
where the terms above are defined by 
\begin{gather}
 B(I) := \widehat{b}(I) {\lambda(I)^{1\over2}  (\mu^{-{1}}(I))^{{1\over2}} \over |I|^{{3\over2}}},
 \\
 H_I(x) := {\mu^{-{1\over2}}(x) |I|^{1\over2} h_I(x)\over  (\mu^{-{1}}(I))^{{1\over2}}},\quad \textup{and} \quad  K_I(x):={\lambda^{1\over 2}(x) |I_+|^{1\over2}k_{I}(x)\over \lambda(I)^{1\over2}}.
\end{gather}
Similar to the estimate in Subsection 3.1, we see that both $\{H_I\}_{I\in\mathcal{D}}$ and $\{K_I\}_{I\in\mathcal{D}}$ are NWOs.

Combining all the estimates above, we see that $[b,\Sh]\in { S^p(L^2_\mu(\R) , L^2_\lambda(\R)) }$. Hence, we further have  $[b,H]\in { S^2(L^2_\mu(\R) , L^2_\lambda(\R)) }$ with
$$\|[b,H]\|_{ S^2(L^2_\mu(\R) , L^2_\lambda(\R)) }\lesssim \|b\|_{B^2_{\nu}(\R)}. $$
The proof of Proposition \ref{prop schattenlarge0} is complete.
\end{proof}
\bigskip

\noindent
{\bf Necessity:}
\begin{proposition}\label{prop schattenlarge1}
Suppose $\mu,\lambda\in A_2$, $\nu= \mu^{1\over2}\lambda^{-{1\over2}}$.
Suppose $b\in {\rm VMO}_\nu(\mathbb R)$ with $\|[b,H]\|_{S^2(L^2_\mu(\R) , L^2_\lambda(\R)) }<\infty$, then we have
\begin{equation}\label{B norm by Sp}
\|b\|_{B_{\nu}^{2}(\R)}\lesssim \|[b,H]\|_{S^2(L^2_\mu(\R) , L^2_\lambda(\R)) }.
\end{equation}
\end{proposition}

\begin{proof}
To begin with, we note that from Lemma \ref{key lemma},
$[b,H]\in S^2(L^2_\mu(\R) , L^2_\lambda(\R))$ if and only if $ \lambda^{1\over 2} [b,H] \mu^{-{1\over2}}$ is in the unweighted $S^2( L^2(\R) , L^2(\R) )$.

Next, we note that from Theorem \ref{thm dyadic structure},  we know that the continuous Besov space is the intersection of two dyadic Besov spaces. 
Thus, it suffices to show that 
$$\|b\|_{B_{\nu,d}^{2}(\R)}\lesssim \|[b,H]\|_{S^2(L^2_\mu(\R) , L^2_\lambda(\R)) },$$
for any dyadic Besov space.  Below, we consider intervals $Q\in \mathcal D$, a fixed dyadic system. 

This argument is of a standard nature. For each dyadic interval $Q$,  let $\hat{Q}$ be the dyadic interval such that:
\begin{enumerate}
  \item $|Q|=|\hat Q|$, and both $Q$ and $\hat Q$ are contained in the same dyadic interval $\widetilde Q$ with
$|\widetilde Q|\lesssim |Q|$. 

\item  The kernel of Hilbert transform $K(x-{\hat x})$ does not change sign for all $(x,{\hat x})\in Q\times\hat{Q}$ and
\begin{equation}\label{lower}
|K(x-{\hat x})|\gtrsim \frac{1}{|Q|}.
\end{equation}

\item Let $\alpha_{\hat{Q}}(b)$ be a median value of $b$ over $\hat{Q}$. This  means $\alpha_{\hat{Q}}(b)$ is a real number such that 
\begin{align} \label{e:E1S}
E_{1}^{Q}:=\{y\in Q:b(y) < \alpha_{\hat{Q}}(b)\}\ \ {\rm and}\ \
E_{2}^{Q}:=\{y\in Q:b(y)>\alpha_{\hat{Q}}(b)\}.
\end{align}
We note that we have  the upper bound $  \lvert  E ^{\hat Q} _{j}\rvert \leq \tfrac{1}2 \lvert  \hat  Q\rvert  $ for $ j=1,2$.
A median value always exists, but  may not be unique. 

\end{enumerate}

Next we decompose $Q$ into two dyadic children by writing $Q=\bigcup_{j=1}^{2}Q_{j}$.
By the cancellation property of $h_{Q}$, we see that
\begin{align}\label{comcom2}
\left|\int_{Q}b(x)h_{Q}(x)dx\right|&=\left|\int_{Q}(b(x)-\alpha_{\hat{Q}}(b))h_{Q}(x)\,dx\right|\nonumber\\
&\leq {1\over |Q|^{1\over2}} \int_{Q}\left|b(x)-\alpha_{\hat{Q}}(b)\right|dx\nonumber\\
&\leq \sum_{j=1}^{2}{1\over |Q|^{1\over2}}\int_{Q_{j}}\left|b(x)-\alpha_{\hat{Q}}(b)\right|dx\nonumber\\
&\leq \frac{1}{|Q|^{1\over2}}\sum_{j=1}^{2}\int_{Q_{j}\cap E_{1}^{Q}}\left|b(x)-\alpha_{\hat{Q}}(b)\right|dx+ \frac{1}{|Q|^{1\over2}}\sum_{j=1}^{2}\int_{Q_{j}\cap E_{2}^{Q}}\left|b(x)-\alpha_{\hat{Q}}(b)\right|dx\nonumber\\
&=:{\rm Term}_{1}^{Q}+{\rm Term}_{2}^{Q}.
\end{align}
Above, we are using the notation \eqref{e:E1S}.

Now we denote
\begin{align*}
F_{1}^{\hat{Q}}:=\{y\in \hat{Q}:b(y)\geq\alpha_{\hat{Q}}(b)\}\ \ {\rm and}\ \
F_{2}^{\hat{Q}}:=\{y\in \hat{Q}:b(y)\leq\alpha_{\hat{Q}}(b)\}.
\end{align*}
Then by the definition of $\alpha_{\hat{Q}}(b)$, we have $|F_{1}^{\hat{Q}}|=|F_{2}^{\hat{Q}}|\sim|\hat{Q}|$ and $F_{1}^{\hat{Q}}\cup F_{2}^{\hat{Q}}=\hat{Q}$. Note that for $s=1,2$, if $x\in E_{s}^{Q}$ and $ y\in F_{s}^{\hat{Q}}$, then
\begin{align*}
\left|b(x)-\alpha_{\hat{Q}}(b)\right|&\leq\left| b(x)-\alpha_{\hat{Q}}(b)\right|+\left|\alpha_{\hat{Q}}(b)-b(y)\right| =\left|b(x)-\alpha_{\hat{Q}}(b)+\alpha_{\hat{Q}}(b)-b(y)\right|= \left|b(x)-b(y)\right|.
\end{align*}
Therefore, for $ s=1,2$, by using \eqref{lower} and by the fact that $|F_{s}^{\hat{Q}}|\approx |Q|$, we have
\begin{align*}
{\rm Term}_{s}^{Q}&\lesssim \frac{1}{|Q|^{1\over2}}\sum_{j=1}^{2}\int_{Q_{j}\cap E_{s}^{Q}}\left|b(x)-\alpha_{\hat{Q}}(b)\right| dx\frac{|F_{s}^{\hat Q}|}{|Q|}\nonumber\\
&= \frac{1}{|Q|^{1\over2}}\sum_{j=1}^{2}\int_{Q_{j} \cap E_{s}^{Q}}\ \int_{F_{s}^{\hat Q}}\left|b(x)-\alpha_{\hat{Q}}(b)\right| \ \frac{1}{|Q|}\ dy\ dx\ \nonumber\\
&\lesssim \frac{1}{|Q|^{1\over2}}\sum_{j=1}^{2}\int_{Q_{j}\cap E_{s}^{Q}}\int_{F_{s}^{\hat Q}}\left|b(x)-\alpha_{\hat{Q}}(b)\right|\left|K(x-y)\right|dydx\nonumber\\
&\lesssim \frac{1}{|Q|^{1\over2}}\sum_{j=1}^{2}\int_{Q_{j}\cap E_{s}^{Q}}\int_{F_{s}^{\hat Q}}\left|b(x)-b(y)\right|\left|K(x-y)\right|dydx.
\end{align*}
To continue, by noting that $K(x-y)$ and $b(x)-b(y)$ do not  change sign for $(x,y)\in (Q_{j}\cap E_{s}^{Q})\times F_{s}^{\hat Q}$, $s=1,2$, we have that
\begin{align}\label{haha}
{\rm Term}_{s}^{Q}&\lesssim \frac{1}{|Q|^{1\over2}}\sum_{j=1}^{2}\left|\int_{Q_{j}\cap E_{s}^{Q}}\int_{F_{s}^{\hat Q}}(b(x)-b(y))K(x-y)dydx\right|\nonumber\\
&=\frac{1}{|Q|^{1\over2}}\sum_{j=1}^{2}\left| \int_{\mathbb R}\int_{\mathbb R}(b(x)-b(y))K(x-y) \chi_{F_{s}^{\hat Q}}(y)\ dy\  \chi_{Q_{j}\cap E_{s}^{Q}}(x)\ dx\right|.
\end{align}
We now insert the weights $\lambda$ and $\mu$ to get
\begin{align*}
{\rm Term}_{s}^{Q}
&\lesssim\frac{1}{|Q|^{1\over2}}\sum_{j=1}^{2}\left| \int_{\mathbb R}\int_{\mathbb R}(b(x)-b(y)) \lambda^{1\over2}(x)K(x-y) \mu^{-{1\over2}}(y) \big(\mu^{{1\over2}}(y)\chi_{F_{s}^{\hat Q}}(y)\big)\ dy\ \big(\lambda^{-{1\over2}}(x)\chi_{Q_{j}\cap E_{s}^{Q}}(x)\big)\ dx\right|. 
\end{align*}

Thus, we further have
\begin{align*}
& \big|\widehat b(Q)\big|^2 {|Q|\over \lambda^{-1}(Q)\ \mu(Q)}  \\
&\lesssim\sum_{j=1}^{2}\left| \int_{\mathbb R}\int_{\mathbb R}(b(x)-b(y)) \lambda^{1\over2}(x)K(x-y) \mu^{-{1\over2}}(y) \big(\mu^{{1\over2}}(y)\chi_{F_{s}^{\hat Q}}(y)\big)\ 
dy\ \frac{ \big(\lambda^{-{1\over2}}(x)\chi_{Q_{j}\cap E_{s}^{Q}}(x)\big) } { \mu(\hat Q)^{1\over2} \lambda^{-1}(Q)^{1\over2}}\ dx\right|^2
\\
&\lesssim \sum_{j=1}^{2}\left| \int_{\mathbb R}\int_{\mathbb R}(b(x)-b(y)) \lambda^{1\over2}(x)K(x-y) \mu^{-{1\over2}}(y) {\mu^{{1\over2}}(y)\chi_{F_{s}^{\hat Q}}(y)\over \mu(\hat Q)^{1\over2}}\ dy\ {\lambda^{-{1\over2}}(x)\chi_{Q_{j}\cap E_{s}^{Q}}(x)\over \lambda^{-1}(Q)^{1\over2}}\ dx \right|^2.
\end{align*}
This implies that
\begin{align*}
\sum_{k}  \sum_{Q\in\mathcal D_{k}} \big|\widehat b(Q)\big|^2 {|Q|\over \lambda^{-1}(Q)\ \mu(Q)} 
&\lesssim  \sum_{k} \sum_{Q\in\mathcal D_k} \sum_{s=1}^{2}\sum_{j=1}^{2}|\big\langle \lambda^{1\over2} [b,H] \mu^{-{1\over2}}G_{\hat Q},
  H_Q \big\rangle|^2.
\end{align*}
On the right above, we are using the notation 
$$ G_{\hat Q}(x):= { \chi_{F_{s}^{\hat Q}}(x) \mu^{{1\over2}}(x)\over \mu(\hat Q)^{1\over2} } \ \ \
{\rm and}\ \ \ H_Q(x)= {\chi_{Q_{j}\cap E_{s}^{Q}}(x)\lambda^{-{1\over2}}(x) \over\lambda^{-1}(Q)^{1\over2} }. 
$$
Both collections of functions abover are NWO for $L^2$.  See Definition \ref{d:NWO}, and Proposition \ref{p:NWO}. 
 And both $Q$ and $\hat Q$ are contained in the same dyadic interval $\tilde Q$ and that the dyadic intervals
$Q, \hat Q$ and $\tilde Q$ have comparable sizes.

It follows from \eqref{e:NWO} that  
$$
\|b\|_{B^2_{\nu,d}(\R)}  \lesssim \|[b,H]\|_{S^2(L^2_\mu(\R) , L^2_\lambda(\R)) }.
$$
The proof of Proposition \ref{prop schattenlarge1} is complete.
\end{proof}

\bigskip
\section{Further discussions on weighted Besov spaces $B^p_\nu(\R)$ for $p\not=2$}
\setcounter{equation}{0}
\label{s: Besov further}

We recall that Peller has characterized the membership of the commutator $[b, H]$ in $S^p$, for all $0< p < \infty $, in the unweighted case. Namely, 

\begin{theorem}\cite{MR1949210}  For a symbol $b\in VMO ( \mathbb R )$, and $0< p < \infty$, we have
\begin{equation}
\lVert [b, H] \rVert_{S ^{p} (L^2 \to L^2 )} ^{p}  \approx  
\int\int \frac{ \lvert  b (x) - b (y) \rvert^p }{ \lvert  x-y \rvert ^2} \; dx dy
\end{equation}
\end{theorem}  		

Our first question is then:  Does this characterization  continue to hold in the one weight case? Namely, when one consider the commuator as acting on $L^2 (w)$, 
for an $A_2$ weight $w$.   Our second question is can our dyadic result Theorem \ref{thm Pi and Haar multiplier} be extended to the case $0 <p \leq 1$?  
Our methods used duality, especially through the inequality \eqref{e:NWO}.  So a new approach would be needed.   

Also, recall that if one consider the commutator with a Riesz transform $R$ in higher dimensions on $\mathbb R ^{d}$, the only symbols $b$ for which the commutator 
is in $S^p$, for $0< p < d$ are the constant symbols. But in that case, the commutator is zero. See the introduction to \cite{MR686178}.  
Thus, given our restriction in our main theorem on the Hilbert transform, the restriction to dimension one is reasonable.  

Our main question is how to extend Theorem \ref{schatten H} to the case of Schatten classes for $p\neq 2$.  Let us describe a difficulty in finding a sufficient condition. 
In our current proof, we used a classical characterization of Hilbert--Schmidt operators in terms of the their kernels. Namely, if $T$, an operator acting on a 
general measure space  has kernel $ K(x,y)$, we have 
\begin{equation}
\lVert T \rVert_{S^2} ^2 = \int\int \lvert  K (x,y) \rvert ^2 \;dx\,dy. 
\end{equation}
This has an extension to $S^p$ for $p>2$.  Combining observations of Russo  \cite{MR500308} and Goffeng \cite{MR2935392}, 
Janson and Wolff \cite{MR686178}*{Lemma 1 and Lemma 2}: if $p>2$ and $1/p+1/p^{\prime}=1$, then
\begin{align}\label{integral}
\|T\|_{S^{p,\infty}}\leq \|K\|_{L^{p},L^{p^{\prime},\infty}}^{1/2}\|K^{*}\|_{L^{p},L^{p^{\prime},\infty}}^{1/2},
\end{align}
where $\|\cdot\|_{L^p, L^{p^{\prime},\infty}}$ denotes the mixed-norm:
$
\|K\|_{L^p,L^{p^{\prime},\infty}}:=\big\|\|K(x,y)\|_{L^p(dx)}\big\|_{L^{p^{\prime},\infty}(dy)}.
$

Consider $T = [b,H]$ and $ \lambda, \mu\in A_2 $.
Then, $T' = M _{ \sqrt{ \lambda}} T M _{1/\sqrt{ \mu}}$
is mapping $L^2 ( \mathbb R)$ to itself with
the kernel 
\begin{equation}
K (x,y) =  \sqrt{\lambda (x)/ \mu (y)} \ \frac {b(x)-b(y)}{x-y} 
\end{equation}
Then, one needs to bound, with $1/q=1-2/p$
\begin{align}\label{verify1}
\lVert K  \rVert_{L ^{p}, L ^{p', \infty}}
& \leq  \Bigl\lVert
\sqrt{\lambda (x)/\mu (y)} \frac{ \lvert b(x)-b(y) \rvert }   { \lvert x-y \rvert  }
 \Bigr\rVert _{L ^{p}, L ^{p', \infty}}
\\
& \leq
\Bigl\lVert
\sqrt{\lambda (x)/\mu (y)} \frac{ \lvert b(x)-b(y) \rvert }   {  |x-y|  } 
 \Bigr\rVert _{L^p, L^p}
 \lVert  \lvert x-y \rvert ^{-1/q} \rVert_{L^ \infty L^ {q, \infty }}
 \nonumber\\
 & \lesssim \Bigl\lVert
\sqrt{\lambda (x)/\mu (y)} \frac{ \lvert b(x)-b(y) \rvert }   { \lvert x-y \rvert ^{2/ p} }
   \Bigr\rVert _{L^p (dxdy)}.  \label{e:possibleBesov}
\end{align}
The expression above is problematic, however.  The power $p/2$ falls on an $A_2$ weight. 
A power weight $ \lambda (x) = \lvert x \rvert ^{- \alpha}$ is an $A_2$ weight for $0< \alpha <1/2$.  For $ \alpha p/2 >1$, this weight will not be in $L^{p/2} _{\textup{loc}}$.   
One might hope that the the expression in \eqref{e:possibleBesov} could be useful, due to the isolated nature of the singularity of the power weight. 
But, we also have this folklore fact: The  singularities of an $A_2$ weight can be dense, as detailed in the Proposition below.  
For $\lambda $ as in Proposition \ref{p:folklore}, if one uses \eqref{e:possibleBesov} as a definition of a Besov space, then only constants are allowed.
This is in sharp contrast to the dyadic Besov spaces, which contain all Haar functions, for instance.

\begin{proposition} \label{p:folklore} 
For any $r >1$, there is a  $\lambda \in A_2$ for which $ \int_I \lambda ^{r}(x) \;dx = \infty $ 
for all intervals $I$. 
\end{proposition}

\begin{proof}
Fix $r>1$  and  set  $ \alpha = (1+1/r)/2$.
Let $0< \delta  < 1/2$, and consider the function 
\begin{equation}
\phi  _{\delta } (x) = 
\begin{cases}
\delta ^{-\alpha },  & n< x < n+ \delta ,\ n\in \mathbb Z 
\\
1,  & \textup{otherwise}
\end{cases}
\end{equation}
This is an $A_2$ weight, with $A_2$ constant independent of $\delta $.  And,  
\begin{equation}
\label{e:BIG}  
\int _{I} \phi _ \delta ^{r} \; dx \geq \delta ^{1-r \alpha } = \delta  ^{-\frac{r-1}{2}} . 
\end{equation}
This becomes unbounded as $\delta \downarrow 0$.

Now, let $n_j = 2 ^{j}$ and $A (1- \alpha ) >2$. Then define 
\begin{equation}
\lambda (x) = \max _{j \geq 1}   \phi _{ 2 ^{-An_j}} (2 ^{n_j} x + 2 ^{-n_j-1}) \eqqcolon \max_{j\geq1} \lambda _j (x) 
\end{equation}
The purpose of having $n_j$ to increase quickly is to have the `peaks' of the weights have rapidly increasing heights, and very short durations. 
The cumulative mass that the peaks add an interval $I$ with length greater than $2^{-n_j}$ is 
\begin{align}  \label{e:SMALL}
2 ^{-An_j (1- \alpha )} 2^{n_j} \lvert  I \rvert \lesssim 2 ^{-n_j} \lvert  I \rvert. 
\end{align}
The shift `$+2^{-n_j-1}$' is to make the peaks disjoint.    
It is a consequence of \eqref{e:BIG} that $\lambda $ has infinite $L^r$ norm on every interval.   We need to check that $\lambda \in A_2$.

Fix an interval $I$.  Certainly, the only interesting case is $\lvert  I \rvert <1$.   
Let $ \mathcal P_j$ be the components of the set $ \lambda _j ^{-1} (2 ^{A \alpha n_j}) $.  Let $ k $ be the minumum integer $j$ 
so that $I$ intersects some interval $P\in \mathcal P_j$. 
If there are no such integers $j$, then $ \lambda $ is identically one on $I$, and there is nothing to prove. 
Otherwise, we have by \eqref{e:SMALL}, 
$
\lambda (I) \lesssim \lambda _k (I) 
$.
And, $ \chi _I  \lambda (x) \geq \chi _I \lambda_k ( x)$, so the $A_2$ property follows from the fact that the $\lambda _k$ 
have uniformly bounded $A_2$ constant.  
\end{proof}

\bigskip
\bigskip
\bigskip
 {\bf Acknowledgements:}
M. Lacey is a 2020 Simons Fellow, his Research is supported in part by grant  from the U.S. National Science Foundation, DMS--1949206. J. Li is supported by the Australian Research Council through the research grant DP220100285.  B. D. Wick's research is supported in part by U. S. National Science Foundation -- DMS 1800057 and Australian Research Council -- DP 190100970.

%\bibliography{Schatten_commutator_Bloom_Hilbert_arXiv2022}  

\end{document}